\input amstex
\documentstyle{amsppt}
\input epsf
\NoBlackBoxes
\magnification 1200
\vcorrection{-9mm}

\topmatter
\title
%
Algebraically unrealizable complex orientations of plane real
pseudoholomorphic curves
\endtitle
\author
       S.~Yu.~Orevkov
\endauthor
\abstract
We prove two inequalities for the complex orientations of a separating
non-singular real algebraic curve
in $\Bbb{RP}^2$ of any odd degree. We also construct a separating non-singular
real (i.e., invariant under the complex conjugation) pseudoholomorphic
curve in $\Bbb{CP}^2$ of any degree congruent to 9 mod 12
which does not satisfy one of these inequalities. Therefore the
oriented isotopy type of the real locus of each of these curves
is algebraically unrealizable. 
\endabstract

\address
IMT, l'universit\'e Paul Sabatier, 118 route de Narbonne, Toulouse, France
\endaddress

\email
orevkov\@math.ups-tlse.fr
\endemail

\address
Steklov Mathematical Institute, Gubkina 8, Moscow, Russia
\endaddress

\address
AGHA Laboratory, Moscow Institute of Physics and Technology, Russia
\endaddress

\endtopmatter

\rightheadtext{Algebraically unrealizable complex orientations}

\def\Z{\Bbb{Z}}
\def\R{\Bbb{R}}
\def\C{\Bbb{C}}
\def\P{\Bbb{P}}
\def\RP{\Bbb{RP}}
\def\CP{\Bbb{CP}}
\def\Lp{\Lambda^{\text{p}}}
\def\Ln{\Lambda^{\text{n}}}

\def\Sep{\operatorname{Sep}}
\def\supp{\operatorname{supp}}

\def\refADK  {1}
\def\refB    {2}
\def\refCop  {3}
\def\refF    {4}
\def\refFO   {5}
\def\refCorr {6}
\def\refG    {7}
\def\refGro  {8}
\def\refH    {9}
\def\refIvSh {10}
\def\refKS   {11}
\def\refM    {12}
\def\refUR   {13}
\def\refOrGAFA {14} \let\refGAFA=\refOrGAFA 
\def\refOrNNGU {15}                         
\def\refOrRE   {16} \let\refAFST=\refOrRE  
\def\refOrAa   {17}                        
\def\refOrWim  {18}
\def\refO      {19} 
\def\refCrel   {20} \let\refOShCrel=\refCrel  
\def\refMMJ    {21} \let\refOShMMJ=\refMMJ    
 \let\refOShAa=\refOSa
\def\refR      {23}
\def\refV      {24}
\def\refWei    {25}
\def\refW      {26}
\def\refWii    {27}

\def\eqSchNine {1}
\def\eqMain    {2}
\def\eqMainBis {3}
\def\eqCongr   {4}
\def\eqPfA     {5}
\def\eqPfB     {6}
\def\eqPfC     {7}

\def\thMain     {1.1}
\def\corNine    {1.2} 
\def\defPH      {1.3}
\def\remAbel    {1.4}
\def\propBigDeg {1.5}
\def\remNineAlg {1.6}
\def\queSwap    {1.7}
\def\remCongr   {1.8}
\def\remSharp   {1.9}
\def\remMinDeg  {1.10}

\def\sectPH    {2}

\def\sectSep   {3}
\def\thG       {\sectSep.1}
\def\thS       {\sectSep.2}
\def\exQuartic {\sectSep.3}
\def\exNine    {\sectSep.4}
\def\exKS      {\sectSep.5}
\def\remSepar  {\sectSep.6}

\def\sectProof {4}

\def\sectConstr {5}            \let\sectBigDeg=\sectConstr
\def\lemBigDeg  {\sectConstr.1}

\def\remBigDeg  {\sectConstr.3}

\def\sectSwap  {6}
\def\propSwap  {\sectSwap.1}

\def\figMain       {1}
 \let\figCubics=\figSchNine
\def\figPerturb    {3}
\def\figNineAlg    {4}
\def\figQ          {5}
\def\figNine       {6}
\def\figKS         {7}
\def\figConstrA    {8}
\def\figConstrB    {9}
\def\figConstrC    {10}
\def\figConstrD    {11}
\def\figSwap       {12}

\document

\head 1. Introduction
\endhead

By a non-singular real algebraic curve in $\RP^2$ we mean a
non-singular algebraic curve in $\CP^2$ invariant under the complex conjugation
$(x:y:z)\mapsto(\bar x:\bar y:\bar z)$. If such a curve is denoted by $A$, then
we denote the set of its real points by $\R A$. A curve $A$ is called
{\it separating} (or {\it type I\/}) if $A\setminus\R A$ is not connected.
In this case $A\setminus\R A$ has two connected components exchanged by
the complex conjugation, and the boundary orientation induced by the
complex orientation of any of these components is called a
{\it complex orientation} of $\R A$. It is defined up to simultaneous
reversing of the orientation of each connected component of $\R A$.

The main result of the paper (Theorem~\thMain\ below)
is an inequality for the isotopy type
of a plane nonsingular real algebraic curve endowed
with a complex orientation (i.e., for the {\it complex scheme} of such curve
according to Rokhlin's terminology [\refR]) which implies in
particular that the oriented isotopy type shown in Figure~\figMain,
that is the complex scheme
(in the notation of Viro [\refV])
$$
   J\sqcup 9_- \sqcup 1_-\langle 1_+\langle 1_-\rangle\rangle       \eqno(\eqSchNine)
$$
is unrealizable by a real algebraic curve of degree $9$ in $\RP^2$.
Since this complex scheme
is easily realizable by a real
pseudoholomorphic curve 
(see Definition~\defPH\ below), it provides the first example
of a complex scheme of a non-singular plane real projective curve
which is algebraically unrealizable but pseudoholomorphically realizable.
We also construct similar examples for any degree congruent to $9$ modulo $12$.

\midinsert
\epsfxsize60mm
\centerline{\epsfbox{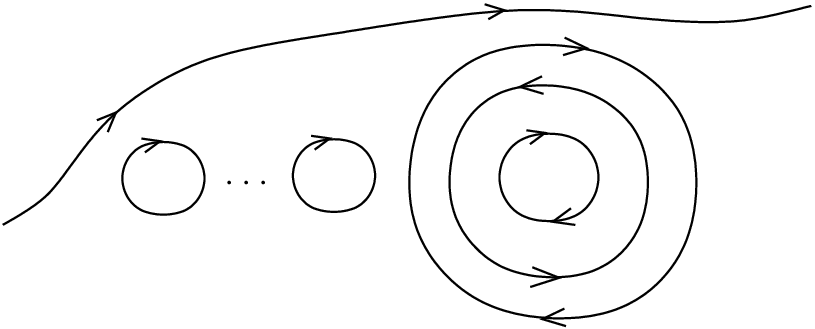}}
\vskip -12mm
\hskip 38mm $\underset{\text{9 negative ovals}}\to{\underbrace{\hskip 17mm}}$
\botcaption{  Figure~\figMain  }
The complex scheme (\eqSchNine)
\endcaption
\endinsert

Let $A$ be a non-singular separating real algebraic curve in $\RP^2$ of an odd degree
$m=2k+1$. We fix a complex orientation on $\R A$.
Let $r$ be the number of connected components of $\R A$. Then $l=r-1$ is
the number of {\it ovals} (components of $\R A$ whose complement in $\RP^2$ is
not connected). The component which is not an oval is called {\it
pseudo-line} and we denote it by $J$.
Following [\refR, \refV], we say that an oval is even (resp. odd) if it is encircled
by an even (resp. odd) number of other ovals.
An oval $O$ is called {\it positive}
if $[O]=-2[J]$ in $H_1(M)$ where $M$ is the closure of the
non-orientable component of $\RP^2\setminus O$. Otherwise $O$ is called {\it negative}.
 Traditionally, the number of
even (resp. odd) ovals is denoted by $p$ (resp. by $n$), and the
number of positive (resp. negative) ovals is denoted by $\Lambda_+$ (resp. $\Lambda_-$).
Let 
$$
\split
       &\Lp_+ = \text{the number of positive even ovals,}\\
       &\Lp_- = \text{the number of negative even ovals,}\\
       &\Ln_+ = \text{the number of positive odd ovals,}\\
       &\Ln_- = \text{the number of negative odd ovals.}
\endsplit
$$

\proclaim{ Theorem \thMain } If $k>0$, then
$$
     \Lp_+ + \Ln_- + 1\ge \frac{l-k^2+2k}2\qquad\text{and}\qquad
     \Ln_+ + \Lp_- \ge \frac{l-k^2+2k}2.                          \eqno(\eqMain)
$$
Setting $l=g-2s$ 
one can equivalently rewrite (\eqMain) in the form
$$
     \Lp_+ + \Ln_- + 1\ge \frac{k^2+k}2-s\qquad\text{and}\qquad
     \Ln_+ + \Lp_- \ge \frac{k^2+k}2-s.                          \eqno(\eqMainBis)
$$
\endproclaim

This theorem is proven in \S\sectProof\ (see  also Example~\exNine).
For the complex scheme (\eqSchNine) we have
$l=12$ and $\Lp_+=\Ln_-=0$, thus the left inequality in (\eqMain) is not satisfied
for $k=4$. So we obtain:

\proclaim{ Corollary \corNine }
The complex scheme (\eqSchNine) is unrealizable by a real algebraic curve
of degree $9$.
\endproclaim

The main interest of Corollary~\corNine\ is that the complex
scheme (\eqSchNine) admits a very simple realization by a real
pseudo-holomorphic curve of degree $9$ which we present just after the following
definition and a brief discussion.

\medskip\noindent
{\bf Definition \defPH.}
Let $(X,\omega)$ be a symplectic $4$-manifold and $c:X\to X$ be
a smooth involution such that $c^*(\omega)=-\omega$.
A {\it real pseudoholomorphic curve} is a $c$-anti-invariant $\Cal J$-holomorphic curve
[\refGro]
for a smooth $c$-anti-invariant almost complex structure $\Cal J$ which is
tamed by $\omega$ (i.e., $\forall\,v\in TX,\,\omega(v,\Cal Jv)>0$).
In this case we denote the fix-point sets of $X$ and $A$ by $\R X$ and $\R A$.
Notice that $\R X$ is a smooth 2-submanifold of $X$ and
$\R A$ is a smooth $1$-submanifold of $\R X$ at smooth points of $A$.
When $X$ is $\CP^2$, we always consider the Fubini-Studi symplectic form
and the standard complex conjugation. In this case we define the {\it degree}
of $A$ as its homological degree, that is $\deg A=m$ if $[A]=m[\CP^1]$ in $H_2(\CP^2)$.

\medskip
In the setting of Definition~\defPH, when a $c$-anti-invariant
$A$ is smooth, it is enough to demand that $A$ is {\it symplectic}
(i.e., $\omega|_A$ does not vanish) because in this case it is necessarily
$\Cal J$-holomorphic for a suitable $\omega$-tame $c$-anti-invariant $\Cal J$.
Indeed, an $\omega$-tame $\Cal J$
is a section of a fibration over $X$ by open balls
[\refGro, Lemma 2.3.${\text{C}}_2$], thus it can be extended from $A$ to $X$
(see [\refWii, Prop.~1.1] for the real case).
Similarly, if $A$ is nodal (each singularity is an intersection
of two smooth transverse local branches), then, in addition to the symplecticity,
it is enough to demand that the intersections are positive
(see [\refIvSh, Lemma~1.4.2] whose proof can be easily adapted for the real case).

\medskip

\midinsert
\epsfxsize120mm
\centerline{\epsfbox{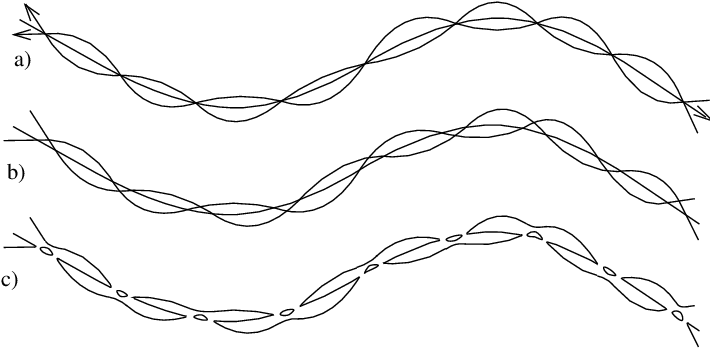}}
\botcaption{  Figure~\figCubics  }
Pseudoholomorphic realization of (\eqSchNine)
\endcaption
\endinsert

Let us show that the complex scheme (\eqSchNine) is realizable by a real
pseudoholomorphic curve in $\CP^2$ of degree $9$ (see Remark~\remNineAlg\
below for another realization).
Let $C=\{f=0\}$ be a real cubic curve with an oval, and $L=\{l=0\}$
be the union of three lines, each
cutting the pseudo-line of $C$ at three distinct real points.
Let $A_{\text{sing}} = \{fg=0\}$ with $g=(f+\varepsilon l)(f-\varepsilon l)$
and $0<\varepsilon\ll 1$.
Then $A_{\text{sing}}$ is a reducible algebraic curve of degree 9 with nine triple points.
Its real locus consists of three nested ovals and a union of three pseudolines arranged
as shown in Figure~\figCubics(a).
In the class of real pseudoholomorphic curves, it can be perturbed as in
Figure~\figCubics(b).
If we consider $f$ and $l$ as
holomorphic sections of the line bundle $\Cal O_{\CP^2}(3)$ rather than
homogeneous polynomials, then the perturbation can be realized by replacing
$f$ with $f+h$ where $h$ is a $\Cal C^1$-small smooth (non-analytic) $c$-invariant
section which is complex analytic in some neighbourhoods of the triple points.
If $h$ is small enough, the obtained curve is analytic near
all double points. Finally, we perturb the double points by adding to $(f+h)g$
a yet smaller $c$-invariant section of $\Cal O_{\CP^2}(9)$
whose signs at the double points are
chosen so that the real locus of the resulting curve $A$ is
the union of three nested ovals with the curve shown in Figure~\figCubics(c).
If the complex orientations of the cubics are chosen as in Figure~\figCubics(a),
the perturbation is coherent with them (see Figure~\figPerturb),
and hence (see  [\refF], [\refM, p.~II.4], or [\refR, \S3.7]) 
the resulting curve of degree $9$ is separating and its complex scheme is
(\eqSchNine).
The non-analytic part of $A$ is close to $A_{\text{sing}}$, hence $A$ is symplectic.

\midinsert
\centerline{\epsfxsize=30mm\epsfbox{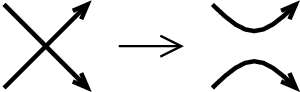}}
\botcaption{Figure \figPerturb}
A perturbation according to complex orientations.
\endcaption
\endinsert

\medskip\noindent
{\bf Remark \remAbel.}
Since (\eqSchNine) is algebraically unrealizable,
so is the intermediate nodal curve. This fact however is much easier:
it immediately follows from Abel's theorem applied to the divisors cut by
any two of the cubic curves on the third one.

\medskip
The same tripling construction applied to suitable Hilbert's $M$-curves yields
similar examples of higher degrees. The following proposition is proven in
\S\sectBigDeg.

\proclaim{ Proposition~\propBigDeg }
For any positive integer $p$ there exists a real pseudoholomorphic
curve of degree $m=12p-3$ with $l=40p^2-38p+10$ ovals such that
$\Lp_+=0$ and $\Ln_-=2p^2-p-1$.
The complex scheme of this curve is unrealizable
by an algebraic curve of the same degree.
\endproclaim

If a real pseudoholomorphic separating curve has two
nested ovals of opposite signs bounding an annulus free of other ovals,
then the complex scheme obtained by reversing the orientations of these
two ovals is also realizable by a separating real pseudoholomorphic curve
of the same degree. I will call this operation {\it swapping of parallel ovals.}
This is a conjugation-equivariant version of Auroux--Donaldson--Katzarkov's
braiding construction [\refADK] (see Proposition~\propSwap\ in
\S\sectSwap\ for more details).

\medskip\noindent
{\bf Remark \remNineAlg.}
The complex schemes
$J\sqcup 9_- \sqcup 1_+\langle 1_-\langle 1_-\rangle\rangle$ and
$J\sqcup 9_- \sqcup 1_-\langle 1_-\langle 1_+\rangle\rangle$
(which are obtained from (\eqSchNine) by swapping the orientations of any two
consecutive nested ovals) are realizable by real algebraic curves of degree $9$.
Indeed, choose the suitable complex orientations on the cubics
in Figure~\figCubics(a), then perturb the cubics generically, and smooth
the double points according to the chosen orientations, see Figure~\figNineAlg.
Thus the swapping of parallel ovals on any of these two algebraic curves
(see Proposition~\propSwap) provides
another pseudoholomorphic realization of the complex scheme (\eqSchNine).
\medskip

\midinsert
\epsfxsize120mm
\centerline{\epsfbox{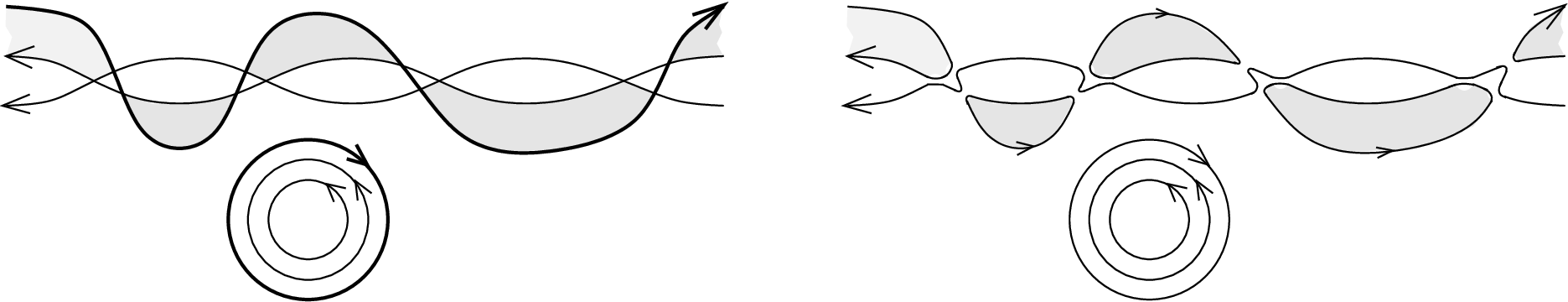}}
\botcaption{  Figure~\figNineAlg  }
Algebraic realization of
$J\sqcup 9_-\sqcup 1_\pm\langle 1_-\langle 1_\mp\rangle\rangle$
\endcaption
\endinsert

\medskip\noindent
{\bf Question \queSwap.}
Is it true that if (\eqMain) does not hold for a separating
real pseudoholomorphic curve $A$, then $A$ can be
transformed to a curve satisfying (\eqMain) by swappings
of parallel ovals?
\medskip

The answer to Question~\queSwap\ is affirmative in all cases considered in this paper
(see Remark~\remNineAlg\ and the proof of Proposition~\propBigDeg).

\medskip\noindent{\bf Remark \remCongr.}
Let $p$, $n$, $\Lambda_+$, and $\Lambda_-$ be the
number of even, odd, positive, and negative ovals respectively.
The arguments in [\refR, \S3.3] (adapted for an odd degree) show that
$p-n\equiv k^2+k+\Lambda_+-\Lambda_-$ mod $4$, whence
$$
  \Lp_+ + \Ln_- \equiv l + \frac{k^2+k}2 \mod 2
   \qquad\text{ and }\qquad
  \Lp_- + \Ln_+ \equiv \frac{k^2+k}2 \mod 2.              \eqno(\eqCongr)
$$
This observation allows us to increase the r.h.s.~of the left (resp.~right)
inequality (\eqMainBis) by $1$ when $l\equiv s$ mod $2$ (resp.~when $s$ is odd).
However, the congruences (\eqCongr) are satisfied by all
separating real pseudoholomorphic curves, thus this improvement
is useless for distinguishing between algebraic and pseudoholomorphic
realizability.

\medskip\noindent
{\bf Remark \remSharp.} (On the sharpness of (\eqMainBis).)
The Rokhlin-Mishachev formula [\refR, \S2.3] implies that the right
inequality (\eqMainBis) is sharp for Harnack curves ($M$-curves
without any nest). One can easily prove by induction that the left inequality
(\eqMainBis) improved according to Remark~\remCongr\ is sharp for Hilbert's
$M$-curves (cf.~\S\sectConstr\ below). I checked that both improved inequalities
are sharp for $k\le 3$ and for any $s$ as soon as this fact does not contradict
the non-negativity of $\Lp_\pm+\Ln_\mp$.

\medskip\noindent
{\bf Remark \remMinDeg.} If $1\le k\le 3$, then all separating real
pseudoholomorphic curves of degree $2k+1$  satisfy the inequalities (\eqMain).
For $k\le2$, this is a consequence of the Rokhlin-Mishachev formula.
For $M$-curves of degree $7$, this is proven in [\refUR, Thm.~2.1],
and the same arguments work for other separating curves of degree $7$.

\medskip
\subhead  Acknowledgements \endsubhead
I am grateful to G.~Mikhalkin,  S.~Nemirovski, and especially to E.~Shustin
for many stimulating discussions. To a great extent this work was inspired by
Kummer and Shaw's paper [\refKS] though I did not use explicitly their results.
I thank the referee for many useful remarks and suggestions.


\head \sectPH. Real algebraic and real pseudoholomorphic curves
\endhead

This section is not used in the rest of the paper and it
can be considered as an extension of the introduction.
A reader interested in the proofs
of the results formulated above can skip it.

It is still unknown if there exists an isotopy type of
configurations of disjoint embedded circles in $\RP^2$
(a real scheme according to Rokhlin's terminology [\refR]) which is
realizable by a smooth pseudoholomorphic curve but algebraically unrealizable with
the same degree.
It seems very plausible that the 6 open cases for $M$-curves of degree 8
(see [\refOrGAFA]) as well as the most of the pseudoholomorphic curves constructed in
[\refOrWim] are such examples but the existing methods are insufficient
to prove it (see the discussion in [\refFO, \S1]).
It is also unknown if there exists a real or complex scheme (again in Rokhlin's sense)
in $\RP^2$ which is
realizable by a flexible curve (in the sense of Viro [\refV]) but
pseudoholomorphically unrealizable.

There are known examples of singular (in particular, reducible)
real pseudoholomorphic curves in $\P^2$ whose real loci are
algebraically unrealizable with the same degree (same degrees of irreducible
components). Some simplest examples are given in [\refFO, \S1].
Notice that the algebraic unrealizability of the Pappus-Ringel arrangement
discussed there, can be deduced from Abel's theorem as in Remark~\remAbel,
if one perturbs one triple of lines into a cubic curve and considers the
divisors which are cut on it by two other triples of lines (of course, a school
geometry proof is also possible).

The following is a brief account of methods used in different settings
to prove the algebraic unrealizability of real pseudoholomorphic curves.

\subhead Hilbert-Rohn-Gudkov method \endsubhead
(See [\refCorr, \refOShCrel, \refOShMMJ, \refOShAa].)
Assuming the existence of an algebraic curve $A$ with a given isotopy type of
$\R A$, one can consider its one-parameter equisingular deformations chosen in such a way
that some quantity (for example the length of a line segment or the area of
some component of the complement to $\R A$) monotonically grows or monotonically
decreases. The monotonicity
ensures that the curve must degenerate. Then one chooses another equisingular
one-parameter family and so on.
This gives a tree of a priori possible degenerations whose leaves
are excluded one-by-one by various topological or algebraic arguments.
Notice that in collaboration with Eugenii Shustin,
by this method we found a very long and
complicated proof of Corollary~\corNine. However the proof was never written, so
we cannot to be sure that it was complete.

\subhead Bezout's theorem for the intersection with unstable curves \endsubhead
(See [\refW].) Isotopy types of some real pseudoholomorphic curves constructed
in [\refW] on real Hirzebruch surfaces are unrealizable algebraically
because their algebraicity contradicts Bezout's theorem for the
number of intersection points with the exceptional curve (a curve $E$ with
$E^2\le -2$) which exists for the integrable complex structure but
does not exist for a generic almost complex structure.

\subhead Trigonal curves \endsubhead (See [\refB, \refCorr].) The construction
from [\refOrRE] provides an algorithm to decide if a given fiberwise isotopy
type is realizable or not by a real algebraic trigonal curve.

\subhead Cubic resolvent of a quadrigonal curve \endsubhead
(See [\refCorr, \refOrNNGU, \refOrAa, \refOShMMJ].)
Using the cubic resolvent, the algebraic unrealizability of a
fiberwise arrangement of a quadrigonal curve can be reduced to that
of a mutual arrangement of a trigonal curve and a line. Then one can try
to prove that the trigonal curve itself is algebraically unrealizable
(as in [\refCorr, \refOrNNGU]) or that its mutual arrangement
with the line is topologically unrealizable (as in [\refOrAa, \refOShMMJ]).

This method can be also applied to curves of higher gonality as follows.
Let us consider a plane real curve as the graph of an $n$-valued function on
the upper half-plane.
Then any four univalent branches of this function can be continued to a
4-valued function on a suitable Riemann surface, and we can study its cubic
resolvent. In this way I proved the algebraic unrealizability in
all cases marked by ``$\not\exists^*$~alg'' in the lists in [\refOrAa]
(unpublished).

\subhead Auxiliary pencils of cubics \endsubhead
(See [\refFO].) A promising idea was to exploit the fact that a pencil of
algebraic cubics through 8 base points on $\P^2$ always has one more
base point, whereas in a family of pseudoholomorphic curves through 8 fixed points,
the 9th crossing point of its members may float. 
However, the implementation of this idea in [\refFO] appeared to be erroneous
(see [\refCorr]) and so far there are no examples where this method allows to
prove the algebraic unrealizability of pseudoholomorphically realizable
isotopy types.


\head \sectSep. Some properties of separating morphisms
\endhead

Let $A$ be an abstract real algebraic curve (a Riemann surface endowed with an
antiholomorphic involution) and $f:A\to\P^1$ a real
(i.e. equivariant under the complex conjugation) morphism.
Following [\refCop, \refKS], we say that $f$ is {\it separating} if $f^{-1}(p)\subset\R A$
for any $p\in\RP^1$. It is clear that if there exists a separating morphism $A\to\P^1$,
then $A$ is separating. The converse is also true and, moreover, the following
estimate takes place (which plays a crucial r\^ole in our proof of Theorem~\thMain).

\proclaim{ Theorem \thG } {\rm(Alexandre Gabard [\refG, Thm.~7.1]).} 
Let $A$ be a smooth connected real algebraic separating curve of genus $g$.
Let $r$ be the number of connected components of $\R A$. Then there exists
a separating morphism $A\to\P^1$
of degree at most $(g+r+1)/2$.
\endproclaim

As shown in [\refCop], the bound $(g+r+1)/2$ is sharp for any fixed $g$ and $r$.
Notice that the restriction to $\R A$ of a separating morphism $A\to\P^1$ of degree $n$
is a covering over $\RP^1$ of degree $n$.
The next theorem is a combination of
the adjunction formula (in terms of Poincar\'e residues)
with the Abel-Jacobi theorem.

\proclaim{ Theorem~\thS }
Let $S$ be a smooth real algebraic surface, $A$ be a smooth irreducible
real separating curve
on $S$, and $D$ be a real divisor on $S$ belonging to the linear system $|A+K_S|$.
Assume that $D$ does not have $A$ as a component. We may always write
$D=2D_0+D_1$ with
a reduced curve $D_1$ and an effective divisor $D_0$.
Let us fix a complex orientation on $\R A$ and
an orientation on $\R S\setminus(\R A\cup\R D_1)$
which changes each time when we cross $\R A\cup\R D_1$ at its smooth point
(this is possible because $D\in|A+K_S|$). The latter orientation induces a
boundary orientation on $\R A\setminus(\R A\cap D_1)$.
Let $f:A\to\CP^1$ be a separating morphism. Then it is impossible that, for some
$p_0\in\RP^1$, the set
$f^{-1}(p_0)\setminus\text{\rm supp}(D)$ is non-empty and
the two orientations coincide at each point of this set.
\endproclaim

\demo{ Proof }
We have $D-A\sim K_S$. So, let $\omega_S$ be a real meromorphic $2$-form
realizing this divisor. Let $\omega$ be the Poincar\'e residue of $\omega_S$ on $A$
(if we have $A=\{F(x,y)=0\}$ and $\omega_S=g(x,y)\,dx\wedge dy$
in some local holomorphic coordinates $(x,y)$ on $S$, then $\omega$ is
the restriction to $A$ of the $1$-form $g\,dx/F'_y\,$).
This is a holomorphic $1$-form on $A$.
By construction, the form $\omega_S$ (after a change of sign if necessary)
defines the orientation on $\R S\setminus(\R A\cup D_1)$ described
in the statement of the theorem,
and $\omega$ defines the induced boundary orientation on
$\R A\setminus(\R A\cap D_1)$.

Let us fix a point $p_0\in\RP^1$, and let $t$ be the coordinate on an
affine chart of $\RP^1$ centered at $p_0$.
Since $f|_{\R A}$ is a covering, the parameter $t$ lifts to a local
parameter near each points of $f^{-1}(p_0)$. This means that, for some
interval containing $0$, there are smooth functions
$p_k:I\to\R A$, $k=1,\dots,n:=\deg(f)$, such that
$f^{-1}(t)=\{p_1(t),\dots,p_n(t)\}$ for any $t\in I$.
Let $v_0=d/dt\in T_{p_0}(\RP^1)$ and let
$v_k=p'_k(0)\in T_{p_k(0)}(\R A)$. Then $v_0=f_*(v_k)$ for each $k=1,\dots,n$.
%
The coordinate $t$ defines a complex orientation on $\RP^1$ which is lifted by $f$
to a complex orientation of $\R A$. Without loss of generality we may assume that
it is the one chosen in the statement of the theorem.

Since all the divisors $f^{-1}(t)$ are equivalent to each other,
the Abel-Jacobi theorem implies that, for any $t\in I$,
$$
   \sum_{k=0}^n \int_{p_k([0,t])} \omega = 0.
$$
Differentiating this identity at $t=0$, we get
$\omega(v_1)+\dots+\omega(v_n)=0$.
It remains to notice that $\omega(v_k)=0$ iff $p_k\in\R A\cap \text{\rm supp}(D)$,
and otherwise the sign of $\omega(p_k)$ is positive (resp. negative) iff
the orientation induced by $\omega$ coincides with (resp. is opposite to)
the chosen complex orientation of $\R A$.
\qed\enddemo

Let us illustrate on simple examples how Theorem~\thS\ works.

\medskip\noindent
{\bf Example \exQuartic. } (Cf.~[\refKS, Examples 2.8 and 3.8].)
In [\refO,~Thm.~1] I proved that {\sl a hyperbolic quartic
curve $A$ in $\RP^2$ does not admit any separating morphism to $\RP^1$
whose restriction to the exterior oval has degree} $1$
(the real locus of a hyperbolic quartic consists of two nested ovals).
This result immediately follows from Theorem~\thS. Indeed, suppose
that such a morphism $f$ exists. Fix $p_0\in\RP^2$ and let $p_1$ be the
only point of $f^{-1}(p_0)$ lying on the exterior oval. Then the elements
of $|A+K|$ are lines. Let $D=D_1$ be a line which separates $p_1$ from the
interior oval (see Figure~\figQ).
The complex orientation of
$\R A$ is shown in Figure~\figQ\ by double arrows. The orientation
defined by $\omega$ (from Theorem~\thS) is shown by ordinary arrows.
This contradicts Theorem~\thS.

\midinsert
\epsfxsize45mm
\centerline{\epsfbox{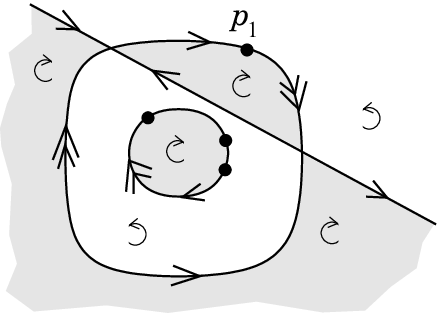}}
\botcaption{  Figure~\figQ  }
Hyperbolic quartic in Example \exQuartic
\endcaption
\endinsert

\medskip\noindent
{\bf Example \exNine.} (A specialization of the proof of Theorem~\thMain\
for the case of the 9th degree complex scheme (\eqSchNine), i.e.,
a direct proof of Corollary~\corNine).
We argue by contradiction.
Suppose (\eqSchNine) is algebraically realizable.
We have $g=28$, $r=13$, $l=12$, $s=8$. 
By Theorem~\thG\ there exists
a separating morphism $f$ of degree $\le(28+13+1)/2=21$.
Elements of the linear system
$|A+K|$ are sextic curves. We choose a double cubic for $D$. Then $D=2D_0$ and $D_1$ 
is empty, thus the orientations are as in Figure~\figNine\ (presented in the
same style as in the previous example).
We see that the two orientations coincide on all ovals of $\R A$ and
they do not coincide on the pseudoline, which we denote by $J$.
Since $l=12$, among the points of $f^{-1}(p_0)$, at least $12$ are on the ovals,
hence at most $9$ of them are on $J$. Thus the cubic curve $D_0$
can be chosen so that $J\cap f^{-1}(p_0)\subset D_0$.
Then we obtain a contradiction with
Theorem~\thS\ unless $f^{-1}(p_0)\subset D_0$. However, the latter case
is impossible. Indeed,
even if there are less than $9$ points in $J\cap f^{-1}(p_0)$,
we can choose any nine-point subset of $J$ containing $J\cap f^{-1}(p_0)$
and trace the cubic $D_0$ through it. Then, if $f^{-1}(p_0)\subset D_0$,
then $D_0$ cuts each oval because each oval has a point of $f^{-1}(p_0)$.
But $D_0$ must cut each oval at an even number of points (counting the multiplicities),
hence it cuts the union of all ovals at least at $24$ points in addition to
the $9$ points where it cuts $J$. This contradicts the Bezout theorem.

\midinsert
\epsfxsize75mm
\centerline{\epsfbox{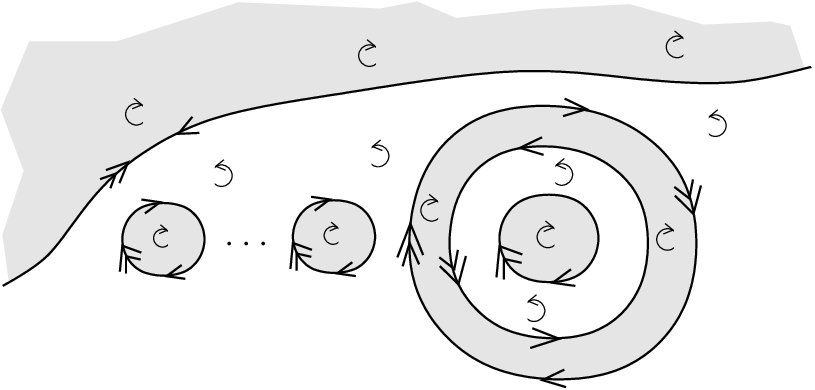}}
\botcaption{  Figure~\figNine  }
The orientations for (\eqSchNine) in Example \exNine
\endcaption
\endinsert

\medskip\noindent
{\bf Example \exKS. } (Cf.~[\refKS, Example 2.14].)
Let $S$ be the complexification of the standard sphere in $\R^3\subset\RP^3$, and
$A$ a curve of degree $6$ on it (a complete intersection of $S$ with
a real cubic surface). Suppose that $\R A$ consists of three nested embedded
circles (see Figure~\figKS), and let $A_0$ be the middle one.
We are going to show that there is no
separating morphism $f:A\to\P^1$ whose restriction to $A_0$
has covering degree $1$. The proof is the same as in Example~\exQuartic.
In this case, the elements of $|A+K|$ are plane sections. So, if such
a morphism exists, we chose
$D=D_1$ to be a circle (cut on $S$ by a plane) which separates the
single point of $A_0\cap f^{-1}(p_0)$ from the other points of $f^{-1}(p_0)$,
and we obtain a contradiction with Theorem~\thS\ (see Figure~\figKS).

\midinsert
\epsfxsize35mm
\centerline{\epsfbox{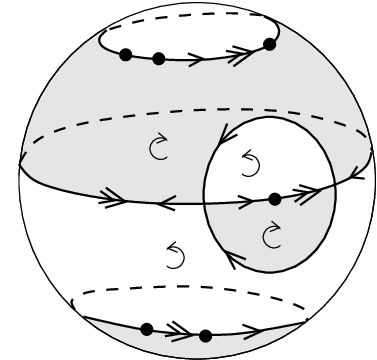}}
\botcaption{  Figure~\figKS  }
The orientations in Example \exKS
\endcaption
\endinsert

\medskip\noindent
{\bf Remark \remSepar.} For an abstract real algebraic curve $A$ with
$r$ components $A_1,\dots,A_r$ of the real locus, Kummer and Shaw
[\refKS] defined the {\it separating semigroup} $\Sep(A)$ of $A$ as the set of
$r$-tuples $(\deg f|_{A_1},\dots,\deg f|_{A_r})$ for all separating morphisms
$f:A\to\P^1$. They proved some properties of it, in particular,
they computed $\Sep(A)$ for all curves $A$ of genus $\le 2$. In [\refO],
I computed the separating semigroup for all hyperelliptic curves and for all curves
of genus $3$. Since each non-hyperelliptic curve of genus $4$ embeds
to a quadric surface $S$ in $\P^3$, Theorem~\thS\ allows to compute $\Sep(A)$ for
any curve $A$ of genus $4$ (cf.~Example~\exKS). In particular,
it appears (though it is not evident a priori) that $\Sep(A)$ for $A$ of genus $4$
depends only on the equivariant deformation class (called also the rigid isotopy type)
of the pair $(S,A)$. These results will be written in a forthcoming paper.

\if{
k0=p
m0=2p+1
m=3m0
k=3p+1
Expand[2k+1-m]
l=3(m0-1)(m0-2)/2+m0^2
Expand[ (l-k^2+2k)/2 ]

}\fi


\head \sectProof. Proof of Theorem~\thMain
\endhead

Let the notation be as in Introduction. So, $A$ is a plane real algebraic non-singular
separating curve of degree $m=2k+1$, and hence of genus $g=(m-1)(m-2)/2=k(2k-1)$.
Let $f:A\to\P^1$ be a separating morphism of degree
$n\le(g+r+1)/2$ which exists by Theorem~\thG.

Suppose that one of the inequalities (\eqMain) does not hold.
Let $C_0$ be the union of the components of $\R A$ which are counted in the
left hand side of that inequality (here we assume that
the ``$1$'' in the left hand side of the left inequality counts the pseudoline of $A$),
and let $C_1=\R A\setminus C_0$. Let $G$ be a homogeneous polynomial
of degree $k-1$, and $D=2D_0$ be the divisor of $G^2$.
Then $\deg D=2k-2=m-3$, i.e., $D\in|A+K_{\P^2}|$. So, let us introduce an orientation
on $\RP^2\setminus\R A$ as in Theorem~\thS\ (notice that $D_1$ is empty in our case).
Then, up to reversing the chosen orientation, we may assume that
the boundary orientation induced from $\RP^2\setminus\R A$ and the complex
orientation are coherent on $C_1$ and not coherent on $C_0$.

For $j=0,1$, let $n_j$ be the covering degree of $f|_{C_j}$ and $r_j$ be the number of
components of $C_j$. Then $r_0$ is the left hand side of the inequality in (\eqMainBis)
which fails, i.~e.,
$$
    r_0\le \frac{k^2+k}2-s-1.                    \eqno(\eqPfA)
$$
We also have $r_1\le n_1=n-n_0$, hence 
$$
   r_0 = r-r_1 \ge r-n + n_0 \ge r+n_0 - \frac{g+r+1}2
       = n_0 + \frac{r-g-1}2
       = n_0 - s.
$$
By combining the two inequalities, we obtain $n_0\le k(k+1)/2 - 1$,
hence we can choose $G$ so that the support of $D$ passes through
$f^{-1}(p_0)\cap C_0$.

Thus, to obtain a contradiction with Theorem~\thS,
it remains to check that the above choice of $G$ can be done so that
$\supp D$ does not pass through
all the $n$ points of $f^{-1}(p_0)$.
Indeed, we may choose $D$ so that
it passes through at least $k(k+1)/2-1$ points of $C_0$.
Suppose that $f^{-1}(p_0)\subset\supp(D)$.
Each component of $C_1$ has at least one point of $f^{-1}(p_0)$, and
at least $r_1-1$ of them are ovals. Since $D_0$ intersects
each oval at least twice, by the Bezout theorem we obtain
$$
    (k-1)(2k+1)\ge D_0\cdot C_0+D_0\cdot C_1 \ge
      \big((k^2+k)/2-1\big) + \big(2r_1-1\big).       \eqno(\eqPfB)
$$
Since $r_0\ge 0$, the inequality (\eqPfA) implies
$$
    s \le (k^2+k)/2-1,          \eqno(\eqPfC)
$$
 hence, denoting $(k^2+k)/2$ by $a$, we obtain
$$
  r_1 =  r-r_0 
      \overset{(\eqPfA)}\to\ge
     (g+1-2s)-(a-s-1)
     = g-a-s+2
      \overset{(\eqPfC)}\to\ge
      g-2a+3.
$$
By plugging this bound for $r_1$ into (\eqPfB), we obtain
$$
    2k^2-k-1 \ge (5k^2-7k+8)/2,
$$
that is $-k^2+5k-10\ge 0$ which is a contradiction.
\if01{
lhs=(k-1)(2k+1); rhs=( a-1 + 2r1-1 )
Factor[rhs//.{r1->g-2a+3,g->k(2k-1),a->k(k+1)/2} ]
Solve[
}\fi
Theorem~\thMain\ is proven.


\head\sectConstr. Construction of curves of degree $12p-3$
(proof of Proposition~\propBigDeg)
\endhead



The recursive construction in the proof of the following lemma
is nothing else than a particular case of Hilbert's construction of $M$-curves in
[\refH] (see also [\refV, \S1.10]), and
we compute the quantities $\Lp_+$ and $\Ln_-$ for the resulting curves.

Let us fix a smooth real conic $E$ on $\P^2$ with $\R E\ne\varnothing$.
Let $\Delta$ be the disk bounded by $\R E$ on $\RP^2$.
For a separating curve $C$, let $\Lambda_{-}(C,\Delta)$
be the number of its negative ovals contained in $\Delta$.

\proclaim{ Lemma~\lemBigDeg }
For any $p\ge 2$ there exists an $M$-curve $C_d$ of degree $d=4p-1$ in
$\RP^2$ such that:
\roster
\item"(i)"
     $C_d$ is transverse to $E$ and
     all the intersection points are real and
     belong to an oval $V$ of $C_d$ which is arranged
     with respect to $\R E$ as shown in Figure~\figConstrA;
\item"(ii)"
     the oval $V$ and all ovals of $C_d\cap\Delta$
     do not encircle other ovals, and the grey
     digons in Figure~\figConstrA\ do not contain ovals of $C_d$;
\item"(iii)"
     $E$ is encircled by $2p-3$ ovals of $C_d$ and $V$ is positive;
\item"(iv)"
     $\Lambda_{-}(C_d,\Delta)=\Ln_-(C_d)=2p^2-p-1$ and $\Lp_{+}(C_d)=0$.
\endroster
\endproclaim

\midinsert
\centerline{\epsfxsize=30mm\epsfbox{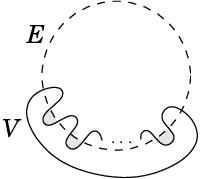}}
\botcaption{Figure \figConstrA} Mutual arrangement of $V$ and $\R E$
\endcaption
\endinsert

\demo{Proof} Induction on $p$. For the base case $p=2$ (i.e., $d=4p-1=7$),
see Figure~\figConstrB.
The inductive step is shown in Figure~\figConstrC.
We see there that when passing from $C_{4p-1}$ to the intermediate
curve $C_{4p+1}$,
the newly appearing ovals are:
\roster
\item"$\bullet$" $4p-1$ positive ovals in $\Delta$,
\item"$\bullet$" $4p-3$ negative ovals outside $\Delta$,
\item"$\bullet$" a negative oval crossing $E$ at $2(4p+1)$ points,
\item"$\bullet$" a positive oval which encircles $E$
      and all the other new ovals.
\endroster
Similarly, when passing from $C_{4p+1}$ to $C_{4p+3}$,
the newly appearing ovals are:
\roster
\item"$\bullet$" $4p+1$ negative ovals in $\Delta$,
\item"$\bullet$" $4p-1$ positive ovals outside $\Delta$,
\item"$\bullet$" a positive oval crossing $E$ at $2(4p+3)$ points,
\item"$\bullet$" a negative oval which encircles $E$
                 and all the other new ovals.
\endroster
Thus, by induction,
$$
   \Lambda_-(C_{4p+3},\Delta) 
   =(2p^2-p-1)+(4p+1) = 2(p+1)^2-(p+1)-1,
$$
and we see that $4p+1$ new ovals contribute to
$\Ln_-$, no new oval contributes to $\Lp_+$,
and the contributions of the old ovals are not changed. 
\qed\enddemo

\midinsert
\centerline{\epsfxsize=124mm\epsfbox{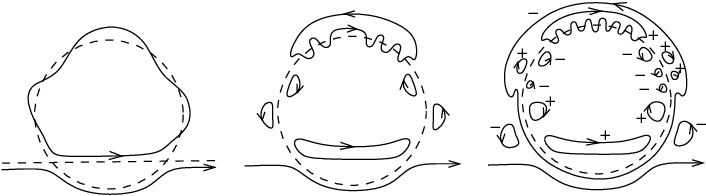}}
\botcaption{Figure \figConstrB} Construction of $C_7$ ($p=2$)
\endcaption
\endinsert

\midinsert
\centerline{\epsfxsize=124mm\epsfbox{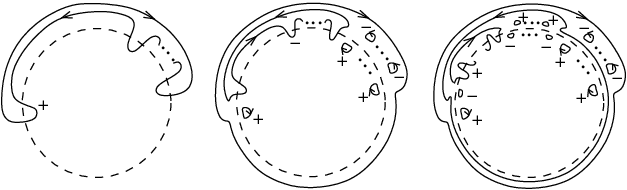}}
\botcaption{Figure \figConstrC} The inductive step in the proof of Lemma~\lemBigDeg
\endcaption
\endinsert

\demo{ Proof of Proposition~\propBigDeg }
We apply the same construction as was explained
in Introduction in the case $p=1$ (see Figure~\figCubics).
Namely, for any $p\ge 2$ and $d=4p-1$, let us consider the curve
$C_d$ from Lemma~\lemBigDeg\ and  let $C'_d$ and $C''_d$ be two its
smooth $\Cal C^1$-small perturbations invariant under the complex conjugation
(by a $\Cal C^1$-small perturbation of a subset $Y$ of a manifold $X$ we mean the
image of $Y$ under a diffeomorphism of $X$ which is $\Cal C^1$-close to identity).
If the preturbations a small enough,  $C'_d$ and $C''_d$ are smooth symplectic
surfaces and we endow them  with the orientation given by the standard symplectic
form on $\CP^2$. Let $A_{\text{sing}}=C_d\cup C'_d\cup C''_d$.
The perturbations can be chosen so that the union of the pseudoline components
of $C_d$, $C'_d$, and $C''_d$ looks as in Figure~\figCubics(a) (with
$d^2$ triple points), the curves are analytic near the triple points,
do not have other intersections,
and for each oval of $\R C_d$ there are ovals of
$\R C'_d$ and $\R C''_d$ which are $\Cal C^1$-close to it.
To acheive these properties, we first construct $C'_d$ as a small
perturbation of $C_d$, and then $C''_d$ as a yet smaller perturbation of $C_d$.
When constructing each of $C'_d$, $C''_d$, we start by
perturbing a neighborhood of $\R C_d$ with the required properties, and then
we extend the perturbation to the whole $C_d$. If crossing points appear
in the complement of $\RP^2$,
the signed number of them is zero, hence they can be removed by pairs by
modifiying the perturbation along paths in $C_d\setminus\R C_d$
connecting crossings of opposite signs.

We choose the complex orientations on the pseudoline components of
$\R A_{\text{sing}}$ as in Figure~\figCubics(a)
and we perturb $A_{\text{sing}}$ to $A$ as in Figure~\figConstrD\
which can be done in two steps similarly to Figure~\figCubics.
Namely, at the first step we perturb only one of the three curves so that
each triple point splits into three double points and so that
the appeared small triangles are placed
on the two sides of any of the three curves in an alternating way.
We may do it so that all the three curves are complex analytic near the double points.
At the second step we perturb the double points
according to the orientations (as in Figure~\figPerturb). This can be done by
a complex analytic perturbation near the double points which is smoothly
extended to the remaining part of the curve respecting the symplecticity.

The perturbation can be chosen so that the ovals of $C'_d$ and $C''_d$ appear
on either side of the corresponding oval of $C_d$, and the side
can be chosen arbitrarily and independently for each oval of $C_d$.
So, we assume that the ovals of $C_d$ contributing to $\Lp_-+\Ln_+$
(let us call them {\it good} ovals) are accompanied by two ovals of $C'_d\cup C''_d$
from the both sides, but for odd positive ovals ({\it bad} ovals),
the both close ovals of $C'_d\cup C''_d$ appear from the interior side
(see Figure~\figConstrD).

\midinsert
\centerline{\epsfxsize=124mm\epsfbox{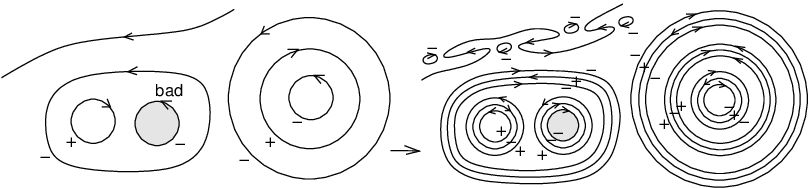}}
\botcaption{Figure \figConstrD} The final step of the construction
\endcaption
\endinsert

Then we have $\Ln_-(A)=\Ln_-(C_d)=2p^2-p-1$ and $\Lp_+=0$ as required.
The number of ovals of $C$ is $d^2$ plus the tripled number of ovals of $C_d$,
that is
$$
    l=3(d-1)(d-2)/2+d^2=40p^2-38p+10
$$
as required.
Then $(l-k^2+2k)/2=2p^2-p+1$ for $k=(m-1)/2=6p-2$,
thus (\eqMain) fails for $A$. Hence the complex
scheme of $A$ is algebraically unrealizable
(we point out that both sides of (\eqMain) are of order $p^2$, but their difference is just $1$).
\qed\enddemo

\medskip\noindent
\newpage

{\bf Remark \remBigDeg.} If we apply the same construction to
the curves of degree $4p+1$ that appear as intermediate curves in the
proof of Lemma~\lemBigDeg, then we obtain the equality sign in the left
inequality in (\eqMain).


\head \sectSwap. Swapping of parallel ovals
\endhead

Let $(X,\omega)$ be a symplectic $4$-manifold and $c:X\to X$ be
a smooth involution such that $c^*(\omega)=-\omega$.
We say that a smooth symplectic surface $A$ in $X$ is {\it real}
if it is anti-invariant under $c$,
i.e., $c(A)=A$ (as sets) and $c_*([A])=-[A]$ in $H_2(X)$.
We denote the fixed point sets of $X$ and $A$ by $\R X$ and $\R A$ respectively.
Then $\R X$ and $\R A$ are smooth submanifolds of $X$
of respective dimensions $2$ and $1$. The condition $c^*(\omega)=-\omega$
implies that $\R X$ is Lagrangian.
We say that $A$ is {\it separating}
if $A\setminus\R A$ is not connected. In this case we define the complex
orientations on $\R A$ in the same way as for real algebraic curves.

The {\it braiding construction} in [\refADK] can be performed
obeying the invariance under $c$ and it provides the following result.
Since the proof in [\refADK] is too sketchy, we give here
a more detailed self-contained proof.

\proclaim{ Proposition \propSwap } {\rm(Essentially, [\refADK, \S3]).}
Let $A$ be a smooth real symplectic surface in $X$.
Let $V_{-1}$ and $V_1$ be two components
of $\R A$ bounding an annulus $B$ on $\R X$. Then there exists
a smooth real symplectic surface $A'$ in $X$ and a neighbourhood $U$ of $B$ in $X$
such that $\R A'=\R A$, $A'\setminus U=A\setminus U$, each of $A\cap U$
and $A'\cap U$ is a union of two annuli
$A_{-1}\cup A_1$ and $A'_{-1}\cup A'_1$ respectively, and
$$
   A_j\cap\R X = A'_{-j}\cap\R X= V_j
   \quad\text{and}\quad A_j\cap\partial U=A'_j\cap\partial U 
 \qquad\text{for $j=\pm1$}.
$$

In particular, if $A$ is separating and $[V_{-1}]=-[V_1]$ in $H_1(B)$
for some complex orientation of $\R A$,
then $A'$ is also separating and
complex orientations on $\R A$ and $\R A'$ can be chosen so that they
are opposite on $V_{-1}\cup V_1$ and
coincide on $\R A\setminus(V_{-1}\cup V_1)$.
\endproclaim


\demo{ Proof }
By Weinstein Neighbourhood Theorem, a symplectic structure
near a Lagrangian submanifold is unique up to symplectomorphism.
Analysing the proof of this result in [\refWei], one can see
that the same is true for $c$-anti-invariant symplectic structures
near the fixed point set of a smooth involution $c$.
Hence we may identify a neighbourhood of $B$ with an open set in
$(\C/\Z)\times\C$ with coordinates
$z_k=x_k+iy_k$, $k=1,2$, where $x_1$ is defined mod $1$,
so that
$$
   B   =\{y_1=y_2=0,\;|x_2|\le 1\},\qquad
   V_j =\{y_1=y_2=0,\; x_2=j\}\qquad (j=\pm1),
$$
$c$ is the usual complex conjugation
$(z_1,z_2)\mapsto(\bar z_1,\bar z_2)$, and $\omega$ is the standard
affine symplectic form $dx_1\wedge dy_1+dx_2\wedge dy_2$.

Let $a,r>0$ be such that the chosen coordinates
are defined in the domain $U_r=\{|y_1|<r, |x_2|<1+a, |y_2|<a\}$.
Since $A$ is $c$-invariant,
for any point $p\in V_j$, the tangent plane $T_p A$ is generated by
the $c$-invariant vector $v_1=\partial/\partial x_1$ 
and a $c$-anti-invariant vector $v_2$ (a linear combination of $\partial/\partial y_1$
and $\partial/\partial y_2$). By combining this fact with $\omega(v_1,v_2)\ne 0$
we conclude that, choosing a smaller $r$ if necessary, we may assume that
the projection $(z_1,z_2)\mapsto z_1$ restricted to
$A_j\cap U_r$ is non-degenerate, that is $A_j\cap U_r$
admits a parametrization
$$
   (u,v)\mapsto\varphi_j(u,v)=\big(u+iv, j+f(u,v) + ig(u,v)\big)
$$
(here $f$ and $g$ depend on $j$).
For any $t<r/2$ we choose a smooth function
$h=h_t:[0,r]\to\R$ such that
$$
     h|_{[0,t]}=0, \quad h|_{[2t,r]}=1,
\quad\text{ and $0\le h'<3/t$ on $[t,2t]$,}
$$
and we modify $A$ in $U_r$ by replacing each $A_j\cap U_r$ with
the annulus $A_{j,t}$ parametrized by 
$$
   (u,v)\mapsto\varphi_{j,t}(u,v)=\big(u+iv,\, j+h(|v|)(f(u,v)+ig(u,v))\big).
$$
We denote the resulting surface by $A_{(t)}$. By construction we have
$c(A_{j,t})=A_{j,t}$ for any $t$.
Let us show that $A_{j,t}$ is symplectic for $t$ small enough.
Indeed, for $v>0$ we have
$$
\xalignat2
   \varphi_j^*(\omega)&=(1+G)du\wedge dv,
     &&G=f'_u g'_v-g'_uf'_v,\\
   \varphi_{j,t}^*(\omega)&=\big(1+h(v)^2G+h(v)h'(v)F\big)du\wedge dv,
                          &&F=f'_ug-g'_uf.
\endxalignat
$$
The smoothness of $f$ and $g$ combined with $f(u,0)=g(u,0)=0$
implies that
$$
     \max_{|v|<2t}\max\big(|f|,|g|,|f'_u|,|g'_u|\big)=O(t).
$$
We also have $|h|<1$ and $|h'|<3/t$
whence $\max_{|v|<2t}|h^2G+hh'F|\to0$ as $t\to0$. Hence we may choose $t$ small enough
such that $A_{j,t}$ is symplectic for $j=\pm1$ and hence so is $A_{(t)}$.
Note that we have $A_{j,t}\cap U_t=\{z_2=j\}$ for $j=\pm1$.

\midinsert
\centerline{\epsfxsize=90mm\epsfbox{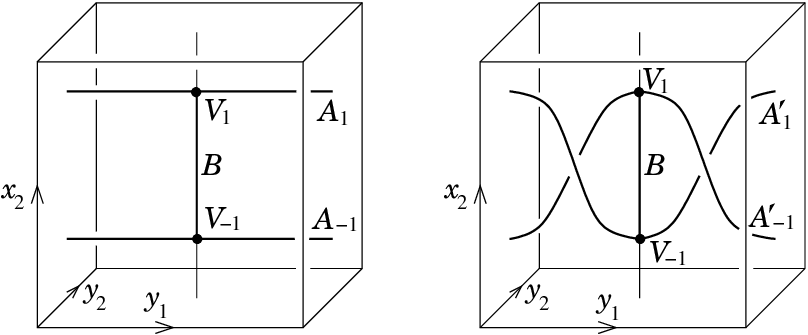}}
\botcaption{Figure \figSwap}
The sections $x_1=\text{const}$ of $A\cup B$ and $A'\cup B$ in $U$
\endcaption
\endinsert

Finally, we define $A'$ as an appropriate smoothing
of the surface obtained from $A_{(t)}$
by replacing $A_{j,t}\cap U_t$ for $j=\pm1$ with the surface
parametrized by
$$
   (u,v)\mapsto\big(u+iv,j\cos(\pi v/t)+jia\sin(\pi v/t)\big)
$$
(see Figure~\figSwap). It is straightforward to check that $A'$ is
symplectic, $c$-invariant, and that it satisfies all other requirements.
\enddemo


\enddocument